\documentclass{article}%
\usepackage{hyperref}
\usepackage{amsmath}
\usepackage{amsfonts}
\usepackage{amssymb}
\usepackage{graphicx}%
\providecommand{\U}[1]{\protect\rule{.1in}{.1in}}
\newtheorem{theorem}{Theorem}

\newtheorem{definition}[theorem]{Definition}

\begin{document}

\title{On the determination of the $grad-div$ criterion}
\author{V. Decaria\thanks{\bigskip The research herein was partially suppported by NSF
grants DMS 1522267 and CBET 160910. vpd7@pitt.edu.}\\Department of Mathematics, \\University of Pittsburgh, \\Pittsburgh, PA 15260, USA
\and W. Layton\thanks{ The research herein was partially suppported by NSF grants
DMS 1522267 and CBET 160910. wjl@pitt.edu, www.math.pitt.edu/\symbol{126}wjl}\\Department of Mathematics, \\University of Pittsburgh, \\Pittsburgh, PA 15260, USA
\and A. Pakzad\thanks{The research herein was partially suppported by NSF grants
DMS 1522267 and CBET 160910. alp145@pitt.edu.}\\Department of Mathematics, \\University of Pittsburgh, \\Pittsburgh, PA 15260, USA
\and Y. Rong\thanks{ rongyao@stu.xjtu.edu.cn. Research of YR supported by NSFC
grants 11171269, 11571274 and Chinese Scholar Council grant 201606280154.}\\School of Math. and Stat.,\\Xi'an Jiaotong Univ.,\\Xi'an, 710049, China
\and N. Sahin\thanks{ nisa70@gmail.com.}\\Department of Mathematics and CS, \\Ankara Yildirim Beyazit Univ., \\Ankara, Turkey
\and H. Zhao\thanks{The research herein was partially suppported by NSF grants DMS
1522267 and CBET 160910. haz50@pitt.edu.}\\Department of Mathematics, \\University of Pittsburgh, \\Pittsburgh, PA 15260, USA}
\maketitle

\begin{abstract}
\textit{Grad-div} stabilization, adding a term $-\gamma\operatorname{grad}%
\operatorname{div}u$, has proven to be a useful tool in the simulation of
incompressible flows. Such a term requires a choice of the coeffiecient
$\gamma$\ and studies have begun appearing with various suggestions for its
value. We give an analysis herein that provides a restricted range of possible
values for the coefficient in $3d$ turbulent flows away from walls. \ If $U,L$
denote the large scale velocity and length respectively and $\kappa$ is the
\textit{signal to noise ratio} of the body force, estimates suggest that
$\gamma$ should be restricted to the range%
\begin{align*}
\frac{\kappa^{2}}{24}LU  &  \leq\gamma\leq\frac{\kappa^{2}}{4}\mathcal{R}%
eLU\text{, mesh independent case},\\
\frac{\kappa^{2}}{24}LU  &  \leq\gamma\leq\frac{\kappa^{2}}{4}\left(  \frac
{h}{L}\right)  ^{-\frac{4}{3}}LU\text{, mesh dependent case.}%
\end{align*}

\end{abstract}
\tableofcontents



\section{Introduction}

\begin{center}
\textit{This is an expanded version of a report with a similar title.}
\end{center}

In the numerical simulation of incompressible flows including an additive term
$-\gamma\operatorname{grad}\operatorname{div}u$ is reported to increase
accuracy \cite{J17}, \cite{LMNOR09}, \cite{O02} and improve performance
of\ iterative solvers \cite{GL89}, \cite{RST08}, \cite{BC99}, \cite{BL12} in
addition to enhancing conservation of mass. In all cases performance depends
on the value of $\gamma$ chosen, \cite{JLMNR16}. Suggested values include:

\begin{itemize}
\item $\gamma=\mathcal{O}(1),$ \cite{BBJL07}, \cite{OLHL09}, \cite{OR04},
\cite{RST08} Remark 3.6 p.308,

\item $\gamma=\mathcal{O}(\nu)$ suggested by estimate (3.17) p.306 in
\cite{RST08},

\item $\gamma=\mathcal{O}(\triangle x^{1\text{ or }2})$, \cite{J17},
\cite{RST08} Remark 3.7 p.308,

\item $\gamma=\mathcal{O}(10^{3})$, \cite{GLRW12}, and

\item global or local ratios of semi-norms $\gamma=\mathcal{O}(|p|_{k}%
/|u|_{k+1})$, \cite{HR13}, \cite{OLHL09}, \cite{J17}, \cite{LAD15}, \cite{H17}.
\end{itemize}

These values of $\gamma$ are derived by optimizing numerical errors for
selected exact solutions, by balancing error estimates for Oseen problems and
by optimizing solver performance. The mesh dependent influence of $\gamma$\ on
total energy dissipation was tested in \cite{CBP16}. We consider herein a
(complementary) mesh independent approach limiting $\gamma$\ to values where
the additional dissipation introduced does not disturb statistical
equilibrium. Denote time averaging by%
\[
\left\langle \phi\right\rangle =\lim\sup_{T\rightarrow\infty}\frac{1}{T}%
\int_{0}^{T}\phi(t)dt.
\]
For $3d,$ fully developed, turbulent flows away from walls, it is known that
the total energy dissipation rate balances energy input, $\left\langle
\varepsilon(u)\right\rangle =\mathcal{O}(U^{3}/L)$\footnote{The energy input
rate at the large scales is $U^{3}/L$. Briefly, the kinetic energy of the
large scales scales with dimensions $U^{2}$. The "\textit{rate}" has
dimensions \textit{1/time}. A large scale quantity with this dimensions is
formed by $U/L$ which is the turn over time for the large eddies, i.e., the
time iit takes a large eddy with velocity $U$ to travel a distance $L$. Thus
the "\textit{rate of energy input}" has dimensions $U^{3}/L$.}, where the
energy dissipation rate (per unit volume) $\varepsilon(u)$ is%
\begin{gather*}
\varepsilon(u)=\frac{1}{|\Omega|}\int_{\Omega}\nu|\nabla u(x,t)|^{2}%
+\gamma|\nabla\cdot u(x,t)|^{2}dx\text{ so}\\
\text{ }\left\langle \varepsilon\right\rangle =\lim\sup_{T\rightarrow\infty
}\frac{1}{T}\int_{0}^{T}\varepsilon(u)dt.
\end{gather*}
This balance is one of the \textit{two laws of experimental turbulence,}
\cite{Frisch} Ch. 5. Building on \cite{DF02}, we analyze the dependence of
$\left\langle \varepsilon\right\rangle $ on $\gamma$ for the \textit{simplest
system arising when incompressibility is relaxed by} $-\gamma
\operatorname{grad}\operatorname{div}u,$ given by%
\begin{equation}
u_{t}+div(u\otimes u)-\frac{1}{2}(\nabla\cdot u)u-\nu\triangle u-\gamma
\nabla\nabla\cdot u=f\,(x)\,.\label{eq:Model}%
\end{equation}
This continuum model arises from the common penalty approximation%
\begin{gather*}
\nabla\cdot u=0\text{ replaced by }\frac{1}{\gamma}p+\nabla\cdot u=0\text{ for
}\gamma>>1\\
\text{so that}\\
\text{ }\nabla p=-\gamma\nabla\nabla\cdot u.
\end{gather*}
Its solutions satisfy the same \'{a} priori energy bound as a discrete NSE
system with $grad-div$ stabilization:%
\[
\frac{1}{2}\frac{d}{dt}||u(t)||^{2}+\left\{  \nu||\nabla u(x,t)||^{2}%
+\gamma||\nabla\cdot u(x,t)||^{2}\right\}  =(f,u).
\]
The domain $\Omega=(0,L_{\Omega})^{3}$ is a $3d$ periodic box, $f(x)$ and
$u(x,0)$ are periodic, satisfy
\[
\nabla\cdot u(x,0)=0,\text{ \ and \ }\nabla\cdot f=0
\]
and have zero mean:%
\begin{gather}
u(x+L_{\Omega}e_{j},t)=u(x,t)\text{ }\\
\text{and}\\
\text{ }\int_{\Omega}\phi dx=0\,\text{\ for }\,\phi=u,\,u_{0},\,f.
\end{gather}
The body force $f(x)$ is assumed smooth so that it inputs energy only into
large scales. Recalling $f(x)$ has mean zero, define the \textit{signal to
noise ratio} of the body force $\kappa$ by%
\[
\kappa=\sqrt{\frac{||f||_{L^{\infty}}^{2}}{\frac{1}{|\Omega|}\int_{\Omega
}|f(x)|^{2}dx}}.
\]
Since $\nabla\cdot u\neq0$ the nonlinearity is explicitly skew symmetrized by
adding $-\frac{1}{2}(\nabla\cdot u)u$.

Let $(\cdot,\cdot),||\cdot||$ denote the $L^{2}(\Omega)$ inner product and
norm. Let $F$, $L$, $U$ denote
\begin{align*}
F &  =\left(  \frac{1}{|\Omega|}||f||^{2}\right)  ^{\frac{1}{2}}\text{,}\\
\text{ }L &  =\min\left\{  L_{\Omega},\frac{F}{||\nabla f||_{L^{\infty}}%
},\frac{F}{(\frac{1}{|\Omega|}||\nabla f||^{2})^{\frac{1}{2}}}\right\}  ,\\
U &  =\left\langle \frac{1}{|\Omega|}||u||^{2}\right\rangle ^{\frac{1}{2}}.
\end{align*}
Non-dimensionalization in the standard way by%
\[
t^{\ast}=\frac{t}{T},\text{ }x^{\ast}=\frac{x}{L},\text{ }U=\frac{L}{T},\text{
}u^{\ast}=\frac{u}{U}%
\]
gives:%
\[
u_{t}^{\ast}+div^{\ast}(u^{\ast}\otimes u^{\ast})-\frac{1}{2}(\nabla^{\ast
}\cdot u^{\ast})u^{\ast}-\frac{\nu}{LU}\triangle^{\ast}u^{\ast}-\frac{\gamma
}{LU}\nabla^{\ast}\nabla^{\ast}\cdot u^{\ast}=\frac{f(x)}{U^{2}}.
\]
We recall $\mathcal{R}e=\frac{LU}{\nu}$ and define the non-dimensional
parameter%
\[
\mathcal{R}_{\gamma}=\frac{LU}{\gamma}.
\]

\begin{theorem}
Let $u(x,t)$ be a weak solution of (\ref{eq:Model}). Then,
\begin{equation}
\left\langle \varepsilon(u)\right\rangle \leq\left(  6+\mathcal{R}e^{-1}%
+\frac{1}{4}\kappa^{2}\mathcal{R}_{\gamma}\right)  \frac{U^{3}}{L}.
\end{equation}

\end{theorem}

This estimate gives insight into $\gamma$ by asking $grad-div$ dissipation be
comparable to (respectively) the pumping rate of energy to small scales by the
nonlinearity, $U^{3}/L$, and to the correction to the asymptotic,
$\mathcal{R}e\rightarrow\infty$, rate due to energy dissipation in the
inertial range, $\mathcal{R}e^{-1}\frac{U^{3}}{L}.$ The cases%
\begin{gather*}
2\simeq\kappa^{2}\mathcal{R}_{\gamma}\\
\text{and}\\
\mathcal{R}e^{-1}\simeq\kappa^{2}\mathcal{R}_{\gamma}%
\end{gather*}
yield
\begin{gather*}
\text{ mesh independent case:\ }\\
\frac{\kappa^{2}}{24}\leq\frac{\gamma}{LU}\leq\frac{\kappa^{2}}{4}%
\mathcal{R}e.
\end{gather*}
Let $\eta\simeq\mathcal{R}e^{-3/4}L$ denote the Kolmogorov miocroscale so
$\mathcal{R}e=(\eta/L)^{-4/3}$. When the model is solved on a spacial mesh
with meshwidth $\eta<<h$ the smallest scale available is $\mathcal{O}(h)$.
Replacing $\eta$ by $h$ leads to an estimate of mesh dependence of%
\begin{gather*}
\text{ mesh dependent case:}\\
\frac{\kappa^{2}}{24}\leq\frac{\gamma}{LU}\leq\frac{\kappa^{2}}{4}\left(
\frac{h}{L}\right)  ^{-\frac{4}{3}}\text{.}%
\end{gather*}

\subsection{Related work}

The energy dissipation rate is a fundamental statistic in experimental and
theoretical studies of turbulence, e.g., Sreenivasan \cite{S84}, Frisch
\cite{Frisch}. \ In 1968, Saffman \cite{S68}, addressing the estimate of
energy dissipation rates, $\left\langle \varepsilon\right\rangle \simeq
U^{3}/L$ , wrote\ that

\begin{center}
\textit{"This result is fundamental to an understanding of turbulence and yet
still lacks theoretical support." }- P.G. Saffman 1968
\end{center}

In 1992 Constantin and Doering \cite{CD92} made a fundamental breakthrough,
establishing a direct link between the phenomenology of energy dissipation and
that predicted for general weak solutions of shear flows directly from the
NSE. This work builds on Busse \cite{B78}, Howard \cite{H72} (and others) and
has developed in many important directions. It has been extended to shear
flows in Childress, Kerswell and Gilbert \cite{CKG01}, Kerswell \cite{K98} and
Wang \cite{Wang97}. For flows driven by body forces extensions include Doering
and Foias \cite{DF02}, Cheskidov, Doering and Petrov \cite{CDP06}\ (fractal
body forces), and \cite{L07} (helicity dissipation). Energy dissipation in
models and regularizations studied in \cite{L02}, \cite{L07}, \cite{LRS10},
\cite{LST10}. 

\section{Analysis of the energy dissipation rate}

Compared to the NSE case \cite{DF02} the term
\[
-\frac{1}{2}(\nabla\cdot u)u
\]
adds dependence on $\gamma$ since $\operatorname{div}u\neq0$. A smooth enough
solution of (\ref{eq:Model}) satisfies the same \'{a} priori energy bound as a
discrete NSE system with $grad-div$ stabilization%
\[
\frac{1}{2}\frac{d}{dt}||u(t)||^{2}+\left\{  \nu||\nabla u(x,t)||^{2}%
+\gamma||\nabla\cdot u(x,t)||^{2}\right\}  =(f,u).
\]
We thus define weak solutions to the model as follows.

\begin{definition}
A weak solution of (\ref{eq:Model}) is a distributional solution satisfying
the energy inequality
\begin{equation}
\frac{1}{2}||u(T)||^{2}+\int_{0}^{T}\nu||\nabla u(t)||^{2}+\gamma||\nabla\cdot
u(t)||^{2}dt\leq\frac{1}{2}||u(0)||^{2}+\int_{0}^{T}(f,u)dt.
\label{eq:ACEnergyEquality}%
\end{equation}

\end{definition}

From (2.1) standard differential inequalities establish that%
\begin{equation}
\frac{1}{2}||u(T)||^{2}+\frac{1}{T}\int_{0}^{T}\varepsilon(u)dt\leq
C<\infty,\text{ }C=C(data)\text{ independent of }T.\label{eq:BoundsU&Epsilon}%
\end{equation}
From (2.2) $\left\langle \varepsilon\right\rangle $\ is well defined and
finite and
\[
\frac{1}{T}||u(T)||^{2}\rightarrow0\text{ as}\ \mathcal{O}(\frac{1}{T}).
\]
$L$ has units of length and satisfies%
\begin{gather}
||\nabla f||_{L^{\infty}}\leq\frac{F}{L}\text{ }\nonumber\\
\text{and}\\
\text{ }\frac{1}{|\Omega|}\int_{\Omega}|\nabla f(x)|^{2}dx\leq\frac{F^{2}%
}{L^{2}}.\nonumber
\end{gather}
Dividing (\ref{eq:ACEnergyEquality}) by $1/(T|\Omega|)$ gives
\begin{gather}
\frac{1}{2T|\Omega|}||u(T)||^{2}+\nonumber\\
+\frac{1}{T|\Omega|}\int_{0}^{T}\nu||\nabla u(t)||^{2}+\gamma||\nabla\cdot
u(t)||^{2}dt\\
\leq\frac{1}{2}\frac{1}{T|\Omega|}||u(0)||^{2}+\frac{1}{T|\Omega|}\int_{0}%
^{T}(f,u)dt.\nonumber
\end{gather}
Define
\[
U_{T}:=(\frac{1}{T}\int_{0}^{T}\frac{1}{|\Omega|}||u||^{2}dt)^{1/2}.
\]
Given (\ref{eq:BoundsU&Epsilon}) and the definition of $F$, this is%
\begin{align}
\frac{1}{T}\int_{0}^{T}\varepsilon(u)dt &  \leq\mathcal{O}(\frac{1}{T}%
)+\frac{1}{T|\Omega|}\int_{0}^{T}(f,u)dt\label{eq:B}\\
&  \leq\mathcal{O}(\frac{1}{T})+F\sqrt{\frac{1}{T}\int_{0}^{T}\frac{1}%
{|\Omega|}||u||^{2}dt}\nonumber\\
&  \leq\mathcal{O}(\frac{1}{T})+FU_{T}\nonumber
\end{align}
To estimate $F$, set the test function in the weak form to be $f(x)$ (recall
$\nabla\cdot f=0$). This yields%
\begin{gather}
F^{2}=\frac{(u(T)-u_{0},f)}{T|\Omega|}-\frac{1}{T|\Omega|}\int_{0}%
^{T}(u\otimes u,\nabla f)-(\frac{1}{2}(\nabla\cdot u)u,f)dt+\label{eq:A}\\
+\frac{1}{T}\int_{0}^{T}\frac{\nu}{|\Omega|}(\nabla u,\nabla f)dt.\nonumber
\end{gather}
Of the four terms on the RHS, by (\ref{eq:BoundsU&Epsilon}) the first
is\ $\mathcal{O}(1/T)$. The second and fourth are bounded using H\"{o}lders
and Young's inequalities by%
\begin{gather*}
\text{second:}\\
\left\vert \frac{1}{T|\Omega|}\int_{0}^{T}(u\otimes u,\nabla f)dt\right\vert
\leq||\nabla f||_{L^{\infty}}\frac{3}{T|\Omega|}\int_{0}^{T}||u||^{2}%
dt\leq3\frac{F}{L}U_{T}^{2},\text{ }\\
\text{fourth:}\\
\left\vert \frac{1}{T}\int_{0}^{T}\frac{\nu}{|\Omega|}(\nabla u,\nabla
f)dt\right\vert \leq\left(  \frac{1}{T}\int_{0}^{T}\frac{\nu}{|\Omega
|}||\nabla u||^{2}dt\right)  ^{\frac{1}{2}}\left(  \frac{1}{T}\int_{0}%
^{T}\frac{\nu}{|\Omega|}||\nabla f||^{2}dt\right)  ^{\frac{1}{2}}\\
\leq\left(  \frac{1}{T}\int_{0}^{T}\frac{\nu}{|\Omega|}||\nabla u||^{2}%
dt\right)  ^{\frac{1}{2}}\frac{\sqrt{\nu}F}{L}\\
\leq\frac{1}{2}\frac{F}{U}\frac{1}{T}\int_{0}^{T}\frac{\nu}{|\Omega|}||\nabla
u||^{2}dt+\frac{1}{2}UF\frac{\nu}{L^{2}}.
\end{gather*}
The third term is treated as follows:%
\begin{gather*}
\text{third:}\\
\left\vert \frac{1}{T}\int_{0}^{T}\frac{1}{|\Omega|}(\frac{1}{2}(\nabla\cdot
u)u,f)dt\right\vert \leq\frac{1}{2}||f||_{L^{\infty}}\sqrt{\frac{1}{T}\int%
_{0}^{T}\frac{1}{|\Omega|}||\nabla\cdot u||^{2}dt}U_{T}\\
\leq\frac{\gamma}{2}\frac{F}{U}\frac{1}{T}\int_{0}^{T}\frac{1}{|\Omega
|}||\nabla\cdot u||^{2}dt+\frac{1}{8\gamma}UF\frac{||f||_{L^{\infty}}^{2}%
}{\frac{1}{|\Omega|}\int_{\Omega}|f(x)|^{2}dx}U_{T}^{2}\\
\leq\frac{\gamma}{2}\frac{F}{U}\frac{1}{T}\int_{0}^{T}\frac{1}{|\Omega
|}||\nabla\cdot u||^{2}dt+\frac{1}{8\gamma}\kappa^{2}UFU_{T}^{2}%
\end{gather*}
Inserting these estimates in (\ref{eq:A}) and simplifying we obtain%
\begin{gather*}
F\leq\mathcal{O}(\frac{1}{T})+\frac{3}{L}U_{T}^{2}+\\
+\frac{1}{2U}\frac{1}{T}\int_{0}^{T}\frac{\nu}{|\Omega|}||\nabla u||^{2}+\\
+\frac{\gamma}{|\Omega|}||\nabla\cdot u||^{2}dt+\\
+\frac{1}{2}\frac{\nu U}{L^{2}}+\frac{1}{8\gamma}U\kappa^{2}U_{T}^{2}.
\end{gather*}
Inserting this in the RHS of (\ref{eq:B}) gives%
\begin{align*}
\frac{1}{T}\int_{0}^{T}\varepsilon(u)dt &  \leq\mathcal{O}(\frac{1}{T}%
)+FU_{T}\\
&  \leq\mathcal{O}(\frac{1}{T})+3\frac{U_{T}^{3}}{L}+\\
&  +\frac{U_{T}}{2U}\frac{1}{T}\int_{0}^{T}\varepsilon(u)dt+\\
&  +\frac{\nu U}{2L^{2}}U_{T}+\frac{1}{8\gamma}U\kappa^{2}U_{T}^{3}.
\end{align*}
Letting $T\rightarrow\infty$ we have, as claimed, that%
\[
\left\langle \varepsilon(u)\right\rangle \leq\left(  6+\mathcal{R}e^{-1}%
+\frac{1}{4}\kappa^{2}\mathcal{R}_{\gamma}\right)  \frac{U^{3}}{L}.
\]

\section{Conclusions}

The analysis herein suggests the following linkage. Weak imposition of
$\nabla\cdot u=0$ at higher Reynolds numbers means explicit skew
symmetrization becomes necessary. Since $\nabla\cdot u\neq0$, this leads to a
second nonlinear term $-\frac{1}{2}(\nabla\cdot u)u$. The parameter $\gamma$
affects the size of $||\nabla\cdot u||$\ which affects the rate at which
$-\frac{1}{2}(\nabla\cdot u)u$ pumps energy to smaller scales. This leads to
restricting $\gamma$ by aligning this energy transfer rate with that of the
underlying incompressible Navier-Stokes equations. To
summarize\footnote{Several helped develop this limeric:
boatertalk.com/forum/LiquidLounge/1052573267.}:

\begin{center}
\textit{Compressibility, however so slight}

\textit{Doubles nonlinearity for skew-symmetry.}

\textit{Cascades can stop by penalty, however light,}

\textit{Unless its criterion is chosen with sagacity.}
\end{center}


\begin{thebibliography}{9999999}                                                                                          %


\bibitem[BIL06]{BIL06}\textsc{L.C. Berselli, T. Iliescu and W. Layton,}
\emph{Large Eddy Simulation}, Springer, Berlin, 2006

\bibitem[B78]{B78}\textsc{F.H. Busse}, \emph{The optimum theory of
turbulence}, Adv. Appl. Mech., 18(1978), 77-121.

\bibitem[BL12]{BL12}\textsc{S. B\"{o}rm and S. Le} \textsc{Borne}, \emph{H-LU
factorization in preconditioners for augmented Lagrangian and grad-div
stabilized saddle point systems.} IJNMF 68(2012) 83--98.

\bibitem[BBJL07]{BBJL07}\textsc{M. Braack, E. Burman, V. John and G. Lube},
\emph{Stabilized finite element methods for the generalized Oseen problem}.
CMAME 196(2007) 853--866.

\bibitem[BC99]{BC99}\textsc{Yu} \textsc{V.\ Bychenkov and E.V. Chizonkov,}
\emph{Optimization of one three-parameter method of solving an algebraic
system of the Stokes type}. Russian J. of Numer. Anal. and Math. Modelling
14.5(1999) 429-440.

\bibitem[CBP16]{CBP16}\textsc{O. Colom\'{e}s, S. Badia and J. Principe,}
\emph{Mixed finite element methods with convection stabilization for the large
eddy simulation of incompressible turbulent flows}, CMAME 304(2016) 294-318.

\bibitem[CDP06]{CDP06}\textsc{A. Cheskidov, C. Doering and N. Petrov},
\emph{Energy dissipation in fractal-forced flow}, J. Math. Phys., 48(2007) 065208.

\bibitem[CKG01]{CKG01}\textsc{S. Childress, R.R. Kerswell and A.D. Gilbert,}
\emph{Bounds on dissipation for Navier-Stokes flows with Kolmogorov forcing},
Phys. D., 158(2001),1-4.

\bibitem[C98]{C98}\textsc{P. Coletti}, \emph{Analytical and numerical results
for k-epsilon and large eddy simulation turbulence models}, Ph.D. Thesis,
UTM-PHDTS 17, U. Trento, 1998.

\bibitem[CD92]{CD92}\textsc{P. Constantin and C. Doering}, \emph{Energy
dissipation in shear driven turbulence}, Phys. Rev. Letters, 69(1992) 1648-1651.

\bibitem[DF02]{DF02}\textsc{C. Doering and C. Foias}, \emph{Energy dissipation
in body-forced turbulence}, J. Fluid Mech., 467(2002) 289-306.

\bibitem[DG95]{DG95}\textsc{C. Doering and J.D. Gibbon}, \emph{Applied
analysis of the Navier-Stokes equations}, Cambridge, 1995.



\bibitem[DKNR12]{DKNR12}\textsc{A.A. Dunca, K.E. Kohler, M. Neda,\ and L.G.
Rebholz}, \emph{A mathematical and physical Study of multiscale deconvolution
models of turbulence}, M2AS, 35(2012) 1205--1219.

\bibitem[F95]{Frisch}\textsc{U. Frisch}, \emph{Turbulence}, Cambridge, 1995.



\bibitem[GLRW12]{GLRW12}\textsc{K. Galvin, A. Linke, L. Rebholz and N.
Wilson}, \emph{Stabilizing poor mass conservation in incompressible flow
problems with large irrotational forcing and application to thermal
convection}, CMAME, 237(2012)166--176.

\bibitem[GL89]{GL89}\textsc{R. Glowinski and P. Le Tallec}, \emph{Augmented
Lagrangian and operator-splitting methods in nonlinear mechanics,} SIAM,
Philadelphia 1989.

\bibitem[H17]{H17}\textsc{N.D. Heavner}, \emph{Locally chosen grad-div
stabilization parameters for finite element discretizations of incompressible
flow problems}, SIURO, 7(2017) SO1278.

\bibitem[HR13]{HR13}\textsc{T. Heister and G. Rapin,} \emph{Efficient
augmented Lagrangian-type preconditioner for the Oseen problem using grad-div
stabilization}, IJNMF, 71(2013)118--134.

\bibitem[H72]{H72}\textsc{L.N. Howard}, \emph{Bounds on flow quantities}, Ann.
Rev. Fluid Mech., 4(1972) 473-494.

\bibitem[J17]{J17}\textsc{E.W. Jenkins, V. John, A. Linke and L. Rebholz},
\emph{On the parameter choice in grad-div stabilization for the Stokes
equations}, Adv. Comp. Math. 40(2014) 491-516.

\bibitem[JK10]{JK10}\textsc{V. John and A. Kindl,} \emph{Numerical studies of
finite element variational multiscale methods for turbulent flow simulations},
CMAME 199(2010): 841-852.

\bibitem[JLMNR16]{JLMNR16}\textsc{V. John, A. Linke, C. Merdon, M. Neilan, and
L. Rebholz,} \emph{On the divergence constraint in mixed finite element
methods for incompressible flows}, SIAM Review (2016).









\bibitem[K98]{K98}\textsc{R.R. Kerswell}, \emph{Unification of variational
methods for turbulent shear flows: the background method of Doering-Constantin
and the mean-flow method of Howard-Busse}, Physica D, 121 (1998), 175-192.

\bibitem[L69]{Lad69}\textsc{O. Ladyzhenskaya}, \emph{The Mathematical Theory
of Viscous Incompressible Flow}, Gordon and Breach, (1969).

\bibitem[L67]{Lad2}\textsc{O. Ladyzhenskaya}, \emph{New equations for the
description of the motions of viscous incompressible fluids, and global
solvability for their boundary value problems. }Trudy Matematicheskogo
Instituta im. VA Steklova 102 (1967): 85-104.

\bibitem[L02]{L02}\textsc{W. Layton}, \emph{Bounds on energy dissipation rates
of large eddies in turbulent shear flows}, Mathematical and Computer Modeling,
35, 2002, 1445 1451.

\bibitem[L07]{L07}\textsc{W. Layton},\emph{\ Bounds on energy and helicity
dissipation rates of approximate deconvolution models of turbulence}, SIAM J
Math. Anal., 39, 916-931 (2007)

\bibitem[LRS10]{LRS10}\textsc{W. Layton, L. Rebholz and M. Sussman},
\emph{Energy and helicity dissipation rates of the NS-alpha and NS-omega
deconvolution models}, IMA Journal of Applied Math, 75, 56-74, 2010.

\bibitem[LST10]{LST10}\textsc{W. Layton, M. Sussman and C. Trenchea,}
\emph{Bounds on energy, magnetic helicity and cross helicity dissipation rates
of approximate deconvolution models of turbulent MHD flows,} Num. Functional
Anal. and Opt., 31, 577-595, 2010.

\bibitem[L54]{L54}\textsc{J. Leray}, \emph{The physical facts and the
differential equations}, American Math. Monthly 61 (1954), 5-7.

\bibitem[LMNOR09]{LMNOR09}\textsc{W. Layton, C. Manica, M. Neda, M.A.
Olshanskii and L. Rebholz}, \emph{On the accuracy of the rotation form in
simulations of the Navier-Stokes equations.} JCP, 228(2009)3433--3447.

\bibitem[LAD15]{LAD15}\textsc{G. Lube, D. Arndt, H. Dallmann,}
\emph{Understanding the limits of inf-sup stable Galerkin-FEM for
incompressible flows}, In: Knobloch P. (eds) BAIL 2014, LN in CSE, vol 10,
Springer, 2015.

\bibitem[MBYL15]{MBYL15}\textsc{W.D. McComb, A. Berera, S.R. Yoffe and M.F.
Linkmann}, \emph{Energy transfer and dissipation in forced isotropic
turbulence}, Phys. Rev. E. 91(2015) 043013.

\bibitem[M97]{Mus96}\textsc{A. Muschinski}, \emph{A similarity theory of
locally homogeneous and isotropic turbulence generated by a Smagorinsky-type
LES}, JFM 325 (1996), 239-260.

\bibitem[O02]{O02}\textsc{M.A. Olshanskii,} \emph{A low order Galerkin finite
element method for the Navier-Stokes equations of steady incompressible flow:
a stabilization issue and iterative methods.} CMAME, 191(2002)5515-5536.

\bibitem[OLHL09]{OLHL09}\textsc{M.A. Olshanskii, G. Lube, T. Heister and J.
L\"{o}we}, \emph{Grad-div stabilization and subgrid pressure models for the
incompressible Navier-Stokes equations}. CMAME, 198(49-52):3975--3988, 2009.

\bibitem[OR04]{OR04}\textsc{M.A. Olshanskii and A. Reusken}, \emph{Grad-Div
stabilization for the Stokes equations}. Math. Comp., 73(2004)1699--1718.

\bibitem[RST08]{RST08}\textsc{H.-G. Roos, M. Stynes and L. Tobiska,}
\emph{Robust numerical methods for singularly perturbed differential
equations: convection-diffusion-reaction and flow problem}s, Springer, Berlin, 2008.

\bibitem[S68]{S68}\textsc{P.G. Saffman}, 485-614 in: \emph{Topics in Nonlinear
Physics}, N. Zabusky (ed.), Springer, 1968.

\bibitem[S01]{Sagaut}\textsc{P. Sagaut,}\emph{Large eddy simulation for
Incompressible flows, }Springer, Berlin, 2001.

\bibitem[S84]{S84}\textsc{K.R. Sreenivasan}, \emph{On the scaling of the
turbulent energy dissipation rate}, Phys. Fluids, 27(1984) 1048-1051.

\bibitem[V15]{V15}\textsc{J.C. Vassilicos}, \emph{Dissipation in turbulent
flows}, Ann. Rev. Fluid Mech. 47 (2015) 95-114.

\bibitem[vD12]{vD}\textsc{E.R. van Driest}, \emph{On turbulent flow near a
wall,} J. of the Aeronautical Sciences (Inst. of Aeronautical Sciences) 23.11 (2012).

\bibitem[vNR50]{vNR}\textsc{J. von Neumann and R.D. Richtmyer}, \emph{A method
for the numerical calculation of hydrodynamic shocks}, J. Applied Physics
21(1950) 232-237.

\bibitem[W97]{Wang97}\textsc{X. Wang, }\emph{The time averaged energy
dissipation rates for shear flows}, Physica D, 99 (1997) 555-563. 2004.
\end{thebibliography}
\end{document}